\documentclass[12pt,draft]{amsart}
\usepackage{amssymb}

\newtheorem{Thm}{Theorem}[section]

\newtheorem{Rem}[Thm]{Remark}
\newtheorem{Exa}[Thm]{Example}

\numberwithin{equation}{section}

\def\al{\alpha}
\def\be{\beta}
\def\Ga{\Gamma}
\def\Om{\Omega}

\def\la{\lambda}
\def\vap{\varphi}
\def\ze{\zeta}

\def\po{\partial\Omega}

\def\C{{\mathbb C}}
\def\D{{\mathbb D}}

\def\R{{\mathbb R}}
\def\wh{\widehat}
\def\wt{\widetilde}
\def\ssm{\smallsetminus}
\def\oli{\overline}
\def\Lolraw{\Longleftrightarrow}
\def\sbs{\subset}
\def\il{\int\limits}

\begin{document}

\title[An Inverse Problem]% for the Double Layer Potential
{An inverse problem for the double layer potential}

\author{P. Ebenfelt$^0$}\footnotetext{The first and second authors were partially supported
by the NSF grants DMS-0100110 and DMS-9703915 respectively}
\address{Department of Mathematics, University of California, San Diego,\newline
La Jolla, CA 92093--0112} \email{pebenfel@math.ucsd.edu}

\author{D. Khavinson$^0$}
\address{Department of Mathematics,
Univeristy of Arkansas,\newline Fayetteville, AR 72701, USA}
\email{dmitry@comp.uark.edu}

\author{H. S.~Shapiro}
\address{Department of Mathematics, Royal Institute of Technology,\newline
S-100~44 Stockholm, Sweden} \email{shapiro@math.kth.se}

%\subjclass{34}

%\keywords{}

\maketitle

\section{Introduction}

\subsection{}The solution to the Dirichlet problem by means of a
double layer potential was initiated by C.~Neumann and H.~Poincar\'e,
and completed in a celebrated paper by I.~Fredholm. What is involved,
very briefly, is the following. (In these introductory remarks we
consider the problem in $\R^n$, although most of the later analysis
will deal with $n=2$. Also, we suppose all data to be as smooth as
needed to justify various assertions.)

If $\Om$ is a smoothly bounded
domain in $\R^n$, $y\in\Om$, and $g(\cdot,y)$ denotes the Green
function of $\Om$ for the Laplace operator with pole at $y$, we have
the well known formula
\begin{equation}
u(y)=\il_{\po} u(x) (\partial/\partial N_x) g (x,y) dS_x
\end{equation}
valid for every harmonic function $u$ in $\Om$ smooth up to the
boundary. Here $\partial/\partial N$ denotes partial
differentiation with respect to the outward directed normal, $dS$
hypersurface measure on $\po$, and $g$ is so normalized that the
{\em Poisson
  kernel}
\begin{equation}
P(x,y):=(\partial/\partial N_x) g(x,y)
\end{equation}
satisfies
\begin{equation}
\il_{\po} P(x,y) dS_x=1.
\end{equation}

If the Poisson kernel for $\Om$ were known, the solution to the
Dirichlet problem
\begin{alignat}{2}
\Delta u &= 0 &&\qquad\textrm{in $\Om$}\\
u &= f &&\qquad\textrm{on $\po$}
\end{alignat}
where, say, $f\in C (\po)$, would be given by
\begin{equation}
u(y) = \il_{\po} f(x) P(x,y) dS_x.
\end{equation}

Now, in general we do not know the Green function (nor the Poisson
kernel), but we know its singular part since
\begin{equation}
g(x,y) = E(y-x) + h(x,y)
\end{equation}
where $E$ is the ``fundamental singularity'' of the Laplace
operator, or ``Newtonian kernel''
\begin{equation}
E(x) = c_n\cdot
\begin{cases}
|x|^{2-n} &\qquad (n\ge 3)\\
\log |x| &\qquad (n=2)
\end{cases}
\end{equation}
(with $c_n$ a normalization factor which depends on the dimension
$n$) and $h(x,y)$ is a harmonic function of $x\in\Om$ for each $y$
(and, in fact, symmetric in $x$ and $y$). Thus, comparing (1.6),
(1.2) and (1.7) it is natural to look for a formula representing
$u$ in the form
\begin{equation}
u(y) = \il_{\po} \vap(x) (\partial/\partial N_x) E(y-x) dS_x
\end{equation}
(a so-called double layer potential), with ``dipole (or doublet)
density'' $\vap$ on $\po$; cf. [K]. Of course we cannot expect
that the harmonic function (1.9) will have the desired boundary
values $f$ on $\po$ if we simply choose $\vap$ as (some constant
times) $f$. Rather, we must introduce an operator $J$ ($z$
denoting a point of $\po$) by
\begin{equation}
(J\vap) (z) :=\lim_{\underset{\scriptstyle{y\in\Om}}{y\to z}} \,
\il_{\partial\Om} \vap(x) (\partial/\partial N_x) E(y-x) dS_x
\end{equation}
and show that $J\vap=f$ is solvable. More precisely (and this is
the essence of Fredholm's solution to the Dirichlet problem): with
suitable regularity hypotheses and choice of Banach space $X$ of
functions on $\po$, $J$ operating on $X$ is equal to $I/2+K$,
where $I$ is the identity and $K$ a compact operator on $X$.
Hence, by Fredholm-Riesz theory its {\em surjectivity} (which
implies the solvability of Dirichlet's problem for data $f$ in
$X$) is a consequence of its {\em
  injectivity}. The latter is relatively easy to check (for details
see [K]). Thus, in the end, this program justifies the intuitive
idea to replace the integral operator with kernel $\partial g
/\partial N$ by that with kernel $\partial E/\partial N$, and a
perturbation argument.

The compact operator $K$ such that $J$ is given by $I/2+K$ can be
represented as an integral operator
\begin{equation}
(K\vap) (z) :=\, \il_{\partial\Om} \vap(x) (\partial/\partial N_x)
E(z-x) dS_x,\qquad z\in\partial\Om.
\end{equation}
The inverse problem referred to in the title concerns the
injectivity of this operator or, equivalently, whether or not
$1/2$ is an eigenvalue of the operator $J$. In two dimensions,
this question turns out to be equivalent to a certain matching
problem for analytic functions. One of our main results, Theorem
3.19, establishes injectivity of $K$ for a special class of
domains in $\R^2$.

For convenience, we shall rescale the double layer potential so
that the kernel of $K$ corresponds to the space of fixed points
rather than the space of eigenfunctions with eigenvalue 1/2.

{\it Acknowledgements.} The authors would like to thank Siv
Sandvik for typesetting this paper. Part of the research was
carried out while the second author was visiting the
Mittag--Leffler Institute. The authors would like to thank the
Mittag--Leffler Institute, and the NSF who partially supported the
visit.

%\bigskip\noindent
\subsection{The two-dimensional case}

To make more precise statements, we now have  to specify our
smoothness assumptions, and the spaces of functions with which we
work, as well as pay attention to normalization constants. We
henceforth confine ourselves to the two-dimensional case and
rescale the operator as mentioned above. The {\em double layer
potential of the dipole density} $F$ (on the boundary $\Ga$ of a
simply connected planar domain $\Om$) is the harmonic function
\begin{equation}
u(z)=\frac1{\pi} \il_\Ga F(\ze) (\partial/\partial N_\ze)
(\log |z-\ze|) ds_\ze
\end{equation}
where $ds$ denotes arclength on $\Ga$. Since our main interest
here does not concern regularity questions we assume henceforth
that $\Ga$ {\em is an analytic Jordan curve} (without cusp
singularities) and, for some fixed $\sigma>0$, $F$ is in {\em the
space $\Lambda_\sigma
  (\Ga)$ of H\"older continuous (of order $\sigma$) complex-valued
  functions on $\Ga$}. Then (1.12) defines a harmonic function $u$ in
$\Om$ (and another one for $z$ in the exterior domain bounded by
$\Ga$).

To solve the Dirichlet problem (for the interior domain) is to
find, for given $f$ in $\Lambda_\sigma (\Ga)$ (say), a
corresponding $F\in\Lambda_\sigma (\Ga)$ such that $u$, defined in
$\Om$ by (1.12), has the boundary values $f$ on $\Ga$. It is easy
to check that for $F$ real-valued,
\begin{equation}
u(z)=\textrm{Re}
\Bigg[\frac1{\pi i} \il_\Ga \frac{F(\ze)}{\ze-z}\,d\ze\Bigg].
\end{equation}

Now, the Cauchy integral
$$
\frac1{2\pi i} \il_\Ga \frac{F(\ze)}{\ze-z}\,d\ze,
$$
where $F$ is any complex valued function in $\Lambda_\sigma
(\Ga)$, defines a pair of analytic functions: one, denoted $f_i$,
is defined and holomorphic in $\Om$ and the other, denoted $f_e$,
in $\Om':=\wh\C\ssm\oli\Om$, where $\wh\C$ denotes the complex
plane (compactified with the point at infinity). We denote by
$A_\sigma (\Om)$ the Banach space of $\sigma$-H\"older continuous
functions on $\oli\Om$, that are analytic in $\Om$; and by
$A_\sigma (\Om')$ the corresponding space on $\Om'$, which satisfy
moreover the requirement that functions in this space vanish at
$\infty$. Then, as is well known from the theory of singular
integrals (see e.g. [D] or [M]),
\begin{equation}
f_i\in A_\sigma (\Om),\qquad f_e\in A_\sigma(\Om').
\end{equation}

For convenience henceforth (since $\sigma$ remains fixed) we denote
simply
\begin{equation}
A:=A_\sigma(\Om),\qquad A':=A_\sigma(\Om').
\end{equation}

We shall also use the notations $A(\Ga)$, $A'(\Ga)$ to denote the
restrictions to $\Ga$ of the spaces $A$, $A'$.

It is well known [M], that
\begin{equation}
f_i(\ze)-f_e(\ze)=F(\ze),\qquad\ze\in\Ga.
\end{equation}

\subsection{The ``Hilbert operator'' $H$}

Following Kerzman and Stein \cite{KS}, we define the {\em Hilbert
operator} $H$ for the domain $\Om$ as the map taking
$F\in\Lambda_\sigma(\Ga)$ to the boundary values of its Cauchy
integral (from inside) on $\Ga$, in other words $HF=f_i|_\Ga$.
Thus, $H$ is a continuous linear map from $\Lambda_\sigma (\Ga)$
to $\Lambda_\sigma(\Ga)$, and it is {\em
  idempotent}: $H^2=H$. Its kernel is $A'(\Ga)$ and its range
$A(\Ga)$.

\begin{Rem} {\rm It is possible, and indeed desirable, to consider
the analogous operator (still denoted $H$) in other spaces than
the H\"older space, e.g. the Hilbert space $L^2(\Ga; ds)$ but then
certain complications may arise, which we will discuss later.}
\end{Rem}

Now, in terms of the Hilbert operator $H$ we can easily represent
the operator corresponding to $J$ in (1.10) on
$\Lambda_\sigma(\Ga)$. Indeed, denoting by $\Pi F$ the boundary
values from $\Om$ of the double layer potential with doublet
density $F$, we have already noted (this is just a rewriting of
(1.13)) that for real-valued $F$:
\begin{align}
\Pi F &=2\textrm{Re}\, HF \\
&=HF + CHF,\notag
\end{align}
where $C$ denotes complex conjugation.

Thus, for $F$, $G$ real we have
\begin{align*}
\Pi (F+i G)
&=
HF+CHF+i(HG+CHG)\\
&=
H(F+i G)+CH(F-iG)\\
&=
(H+CHC) (F+i G)
\end{align*}
so we have
\begin{align}
\Pi
&=
H+CHC\\
&= H+\wt{H},\notag
\end{align}
where
\begin{equation}
\wt{H}:=CHC
\end{equation}
is an idempotent operator that projects $\Lambda_\sigma(\Ga)$ on
its subspace $CA(\Ga)$.

We can now state the question which is the principal concern of this
paper:

{\em For which choices of $\Ga$ does the operator $\Pi$} (which of
course depends on $\Ga$ as do $H$, etc. although we suppress this in
the notation) {\em admit a non-trivial fixed point?}

Note that $\Pi F=F$ means that, with the special choice of
boundary data $F$, the double layer potential with this same
density gives the solution to Dirichlet's problem. For instance,
if $\Ga$ is a circle, then for all $F$ in a space of codimension
one, indeed for any $F$ with mean value zero, (1.12) solves the
Dirichlet problem for boundary values $F$. This is all well known,
and easy to check. Now, we have the following:

\medskip\noindent
{\bf Theorem 1.20.} {\em $\Pi$ admits a non-trivial fixed point if
  and only if there exist non-constant $f\in A(\Om)$, $g\in A(\Om')$
  with $f=\bar g$ on $\Ga$. Moreover, if there exist such $f$ and $g$, then
  both $f$ and $g$ themselves are
  fixed points of $\Pi$.
}

\begin{proof} (i) Suppose $F$ is a non-constant fixed point. Then
\begin{equation}
HF+\wt{H}F=F.\tag{1.21}
\end{equation}
We consider now two cases:

(a) $\wt{H}F=0$. In this case $H\oli{F}=0$, so $\oli{F}\in
A'(\Ga)$. Also, from (1.21) we have $HF=F$, showing $F\in A(\Ga)$.
Thus $A(\Ga)$ and $C A'(\Ga)$ have a non-trivial common element,
as was to be shown.

(b) $\wt{H}F\not=0$. In this case, applying $H$ to (1.21) gives
$H\wt{H}F=0$. Hence $\wt{H}F\in A'(\Ga)$, and $\wt{H}F$ is not
constant. On the other hand, $\wt{H}F=CH\oli{F}$ is the complex
conjugate of an element of $A(\Ga)$, so we conclude that
$A'(\Ga)\cap CA(\Ga)$ contains a non-constant function. This
concludes the proof of the ``only if'' assertion in Theorem 1.20.

\smallskip

(ii) In the other direction, suppose $f\in A(\Ga)$, $g\in A'(\Ga)$
and $f=\bar g$. We shall show $\Pi f=f$ and $\Pi g=g$. Indeed,
$$
\Pi f=(H+CHC) f = Hf+CH g=f,
$$
and $$ \Pi g=(H+CHC)g=Hg+CHf=g,$$ concluding the proof.
\end{proof}

\medskip

We remark that the proof of Theorem 1.20 above also shows that if
$F$ is a non-trivial fixed point of $\Pi$, then either (a)
$H\oli{F}=0$, in which case $F\in A(\Ga)$ and $\oli{F}\in
A'(\Ga)$, or (b) $H\oli{F}$, which belongs to $A(\Ga)$, is also a
non-trivial fixed point and $CH\oli{F}=\wt{H}F$ is in $A'(\Ga)$.

Observe that the set of fixed points of $\Pi$ is a linear space.
Theorem 1.20 has the following remarkable consequence (which was
also proved, using different methods, in \cite{IK}).

\medskip\noindent
{\bf Corollary 1.22.} {\em If the space of fixed points of $\Pi$
contains a non-constant element, it is infinite dimensional.}

\begin{proof} We shall make use of the elementary fact that
  $\Lambda_\sigma(\Ga)$ is a subalgebra of $C(\Ga)$. Suppose $F$ is a
  non-constant fixed point of $\Pi$. According to the proof of
  Theorem 1.20 (and the subsequent remark) either
\begin{alignat*}{2}
&\textrm{(a)}
&&\quad F\in A(\Ga)\qquad\textrm{and}\qquad\oli{F}\in A'(\Ga),\\
\intertext{or} &\textrm{(b)} &&\quad CH\oli{F} \textrm{\,\, is in
$A'(\Ga)$ and non-constant.}
\end{alignat*}

In case (a), all the powers $F, F^2,  F^3,\dots$ are in $A(\Ga)$
and linearly independent over $\C$, while the corresponding
conjugates $\oli{F}, \oli{F}^2, \dots$ are in $A'(\Ga)$. Therefore
(by Theorem 1.20) $F, F^2, F^3, \dots$ (as well as their complex
conjugates) all are fixed points of $\Pi$.

In case (b), the functions $(CH\oli{F})^n$ $(n=1, 2, \dots)$ are
all in $A'(\Ga)$ and linearly independent. Moreover, each of these
is in the range of $CH$, i.e. is the complex conjugate of an
element of $A(\Ga)$, namely $(H\oli{F})^{n}$. Hence all
$(H\oli{F})^{n}$ $(n=1, 2, \dots)$ (as well as their complex
conjugates) are fixed points of $\Pi$.

Summing up, we have shown that if $F$ is a non-trivial fixed point
of $\Pi$, then either all powers $F^n$ $(n=1, 2, \dots)$, or
$(H\oli{F})^n$ ($n=1,2,\ldots$) are also non-trivial fixed points
of $\Pi$ . This implies the corollary.
\end{proof}

Thus, the problem of fixed points of $\Pi$ is equivalent to a
purely function-theoretic one which we henceforth shall call ``{\em
  the matching problem}'': to find non-constant $f\in A(\Om)$, $g\in
A'(\Om)$ which ``match'' on $\Ga$ in the sense that they are complex
conjugates of one another there.

\begin{Rem} {\rm That $f$ is a fixed point of $\Pi$ is the same as saying,
$f$ is an eigenvector of $\Pi$ corresponding to the eigenvalue
$\la=1$. Thus, whenever 1 is an eigenvalue it is of infinite
multiplicity! No other value of $\la$ can be an eigenvalue of
$\Pi$ of infinite multiplicity because, with some regularity of
$\Ga$, as shown by Fredholm, $\Pi=I+2K$ with $K$ compact and
$K\not=0$ (given by the integral operator (1.11)), so the only
possibility for an eigenvalue of $K$ to be of infinite
multiplicity is if it is 0.

We remark also that, apart from the case where $\Om$ is a disk
already noted, the set of fixed points of $\Pi$ never has finite
codimension. This was shown in [S, Theorem 7.6].}
\end{Rem}

Henceforth, we will no longer explicitly refer to potentials nor the
operator $\Pi$, just the matching problem.

Moreover, from a function-theoretic point of view, the H\"older
continuity assumptions, so useful in potential theory, are
somewhat unnatural, so we shall henceforth consider our question
in the following form:

\medskip\noindent
{\bf Matching problem.} Determine for which Jordan domains $\Om$ there
exist non-constant functions $f$, $g$ such that
\begin{enumerate}
\item[(i)]{$f$ is holomorphic in $\Om$ and extends continously to
    $\Ga=\po$;}
\smallskip
\item[(ii)]{$g$ is holomorphic in $\Om':=\wh\C\ssm\oli\Om$,
    $g(\infty)=0$ and $g$ extends continuously to $\Ga$;}
\smallskip
\item[(iii)]{$f(\ze)=\oli{g(\ze)}$,\qquad $\ze\in\Ga$.}
\end{enumerate}

\begin{Rem} {\rm We shall often impose further restrictions, e.g. that $\Ga$
is {\em analytic}. A very interesting question which, however, we
shall not consider in this paper, is: {\em How much (if any)
  regularity of $\Ga$ is forced by the existence of a non-trivial
  matching pair?}}
  \end{Rem}

To appreciate the role played by regularity hypotheses, observe that
the proof of Corollary 1.22 required that the underlying space of
functions $F$ on $\Ga$ be stable under multiplication. Now, the
notions of double layer potential, Hilbert operator, etc. extend
perfectly well to the setting of densities $F$ that need not be
continuous on $\Ga$, but merely (say) in $L^2(\Ga; ds)$ (and the
corresponding space of holomorphic functions are then the {\em Hardy
  spaces} $H^2(\Om)$, $H^2(\Om')$). But since $L^2$ is not an algebra
the proof we gave that, if 1 is an eigenvalue of $\Pi$ it has
infinite multiplicity, no longer works in this setting, even though
the {\em assertion} remains meaningful. We do not know if it is true.

\section{Some domains for which the matching problem has non-trivial
  solutions}

Thus far, the only known domains allowing non-trivial solutions
for the matching problem are {\em lemniscates}. These lemniscate
solutions (more precisely, the following theorem) were discovered
by Mark Melnikov, who communicated the result to one of the
present authors, and generously consented to its inclusion in our
paper.

\begin{Thm} {\em(M.~Melnikov, unpublished)} Let $R(z)$ be a rational
  function such that $\Ga:=\{z\in\C:|R(z)|=c\}$, where $c>0$, is a
  Jordan curve. Suppose moreover that
\begin{enumerate}
\item[(i)]{$R$ has no poles in $\Om$ (the interior domain bounded by
    $\Ga$)};
\smallskip
\item[(ii)]{$R$ has no zeroes in $\Om'$ (the exterior domain)};
\smallskip
\item[(iii)]{$R(\infty)=\infty$.}
\end{enumerate}
Then, the matching problem for $\Ga$ has non-trivial solutions.
\end{Thm}

\begin{proof}

The pair $R$, $c^2/R$ satisfy the matching requirements.
\end{proof}

We call these the {\em lemniscate solutions} to the matching
problem. We do not know any other solutions.

\begin{Exa} {\rm Let $P$ be a polynomial of degree $\ge1$. Then, for
sufficiently large $c$, $\Ga:=\{z\in\C:|P(z)|=c\}$ is a Jordan
curve containing all zeroes of $P$. Hence the pair $P$, $c^2/P$
solve the matching problem for $\Ga$. (It is easy to construct
also other, non-polynomial, lemniscate solutions).}
\end{Exa}

\section{Some domains for which the matching problem admits only the
  trivial solution}

\subsection{Preliminary material}

\subsubsection{The Schwarz function} Our construction of domains that do not allow
non-trivial solutions of the matching problem requires some rather
delicate preliminary results concerning {\em anticonformal
reflection} with respect to algebraic curves. These results are
all known from earlier work of Avci, Ebenfelt, Gustafsson etc.
(These, and other relevant references, may be found in [S]). Since
they are often embedded in more general considerations concerning
so-called quadrature domains, and afflicted with complications due
to multiple connectivity, etc. we give here a simple, and mostly
self-contained account of the results we shall need. First, let
\begin{equation}
\Ga=\big\{(x,y)\in\R^2:P(x,y)=0\big\},
\end{equation}
where $P$ is a polynomial with real coefficients, which we also
express by $P\in\R[x,y]$, assumed irreducible in the ring
$\C[x,y]$. If we write $x=(z+\bar z)/2$, $y=(z-\bar z)/2i$, the
equation in (3.1) takes the form
\begin{equation}
Q(z,\bar z)=0, \end{equation}where \begin{equation}
Q(z,w):=P\big((z+w)/2, \, (z-w)/2i\big).
\end{equation}
Thus, $Q\in\C[z,w]$ and it is easy to see that $Q$ is irreducible
in this ring. Thus
$$
V:=\big\{(z,w)\in\C^2:Q(z,w)=0\big\}
$$
is an irreducible one-dimensional complex algebraic variety, and
$\Ga$ can be identified as the intersection of $V$ with the ``real
plane'' $\{w=\bar z\}$.

In a neighborhood of a non-singular point $(x_0, y_0)$ of $\Ga$ we
can (writing $x+iy=z$, \, $x_0+iy_0=z_0$) represent $\Ga$ by the
equation $\bar z=S(z)$ where $S$ is the {\em Schwarz function} of
$\Ga$ (cf. [S]). From (3.2)
\begin{equation}
Q\big(z, S(z)\big)=0.
\end{equation}
This holds initially for complex $z$ near $z_0$, and by analytic
continuation, identically. Thus, $S$ is an algebraic function of $z$.

Let us write
\begin{equation}
Q(z,w)=a_0(z)+a_1(z)w+\dots+a_n(z)w^n,
\end{equation}
where the polynomial $a_n(z)$ does not vanish identically. Considered
as a polynomial in $w$, (3.5) has a discriminant $D(z)$, a polynomial
whose zeroes give all possible branch points of $S$ (i.e. a finite
subset of $\C$ which contains all points ``above'' which, in one or
more sheets of its Riemann surface, $S$ has a branch point). Let us
denote the totality of these last-mentioned points by $B$. Then, to
each $z\in\C\ssm B$ there are precisely $n$ distinct values which $S$
takes at $z$, upon analytic continuation to $z$ along a suitable
path. (One of these values may be $\infty$, in case $a_n(z)=0$).

It is often geometrically convenient to work with the antianalytic
function
\begin{equation}
R(z):=\oli{S(z)}
\end{equation}
which we call the (multi-valued) {\em anticonformal reflection}
(ACR) associated with $\Ga$.

For $z\in\C\ssm B$ we denote by $\{S(z)\}$ the unordered set (with
$n$ elements) comprising the values $S$ takes on at $z$, and
likewise for the symbol $\{R(z)\}$. We shall need the following
known result. For the convenience of the reader, we supply a
proof.

\medskip\noindent
{\bf Reciprocity Theorem.} {\em If $z_1, z_2\in\C\ssm B$, then
  $z_1\in\{R(z_2)\}$ if and only if $z_2\in\{R(z_1)\}$.}

\begin{proof}

For any $z,w$ in $\C\ssm B$,
$$
w\in\big\{R(z)\big\}\Lolraw \bar w \in \big\{S(z)\big\}\Lolraw
Q(z,\bar w)=0 \qquad\textrm{(by (3.4))}
$$
and by the same token,
$$
z\in \big\{R(w)\big\}\Lolraw Q(w,\oli z)=0.
$$
So, the result follows from the known fact that,

\begin{equation}
Q(z,\bar w)=\oli{Q(w,\bar z)},\qquad z,w\in\C\end{equation} or,
what is equivalent,
\begin{equation}Q(z,w)=\oli{Q(\bar w, \bar z)},\qquad z,w\in\C,
\end{equation}
which reflects the fact that $Q(z,w)$ is constructed via (3.3)
from a polynomial $P(x,y)$ with real coefficients.

To prove (3.8), simply note that $Q(z,\bar z)=P(x,y)$, where
$z=x+iy$, and, in particular, $Q(z,\bar z)$ is real, i.e.
\begin{equation}
Q(z,\bar z)=\oli{Q(z,\bar z)},\qquad z\in \C.
\end{equation}
If we write $Q^\#(z,w)$ for the (holomorphic) polynomial
$\oli{Q(\bar z,\bar w)}$, then (3.9) can be written
\begin{equation}
Q(z,\bar z)=Q^\#(\bar z,z),\qquad z\in \C,
\end{equation}
which implies, since both sides are analytic functions of $z$ and
$\bar z$, that \begin{equation}Q(z,w)=Q^\#(w,z),\qquad
z,w\in\C.\end{equation} The latter is equivalent to (3.8).

\
\end{proof}

%Recall that we denote by $B$ the (finite) set of points in $\C$ over
%which (in some sheet of its Riemann surface) $S$ has a branch
%point. It is important in the following to show

%
%\medskip\noindent
%{\bf Lemma 3.11.} {\em If $z_0=x_0+iy_0\in K\cap B$, then $z_0$ is a
%  singular point of $K$, i.e. $\partial P/\partial x$ and $\partial
%  P/\partial y$ vanish at $(x_0,y_0)$, where $P$ is the polynomial in
%  (3.1).}

%\begin{proof} $z_0$ is in $B$ if and only if $Q(z_0,w)$ (where $Q$ is
%  defined by (3.3)), as a polynomial in $w$, has a multiple root. This
%  is equivalent to: There exists $w_0\in\C$ such that
%$$
%Q(z_0,w_0)=0\qquad,\qquad(\partial/\partial w) Q (z_0,w_0)=0.
%$$
%>From (3.3) this is equivalent to the vanishing of $P$ and $(\partial
%P/\partial x)+i(\partial P/\partial y)$ for $x=x_0$, $y=y_0$ where
%$x_0=(z_0+w_0)/2$, $y=(z_0-w_0)/2i$. If moreover $z_0\in K$ then
%$w_0=\bar z_0$ so $(x_0,y_0)$ is a point in $\R^2$, hence $\partial
%P/\partial x$ and $\partial P/\partial y$ each vanish at this point.
%\end{proof}

\subsubsection{Quadrature domains}

For our purposes in this paper, we shall only require some properties
of a simple subclass of what are called ``quadrature domains''. For
this subclass, which we call here $R${\em -do\-mains}, the properties
needed can be developed very simply without outside references, and we
shall do so in this section.

\medskip\noindent
{\bf Definition.} {\em An $R$-domain is a simply connected plane
domain which is the image of the unit disk $\D$ under the
conformal mapping by a rational function which is univalent in a
neighborhood of $\overline{\D}$ .}

\medskip

Note that our definition implies that the boundary $\Ga$ of
$\Om:=\vap(\D)$ (where $\vap$ is the rational conformal map) is an
algebraic Jordan curve without singular points. What is of interest
for us here is certain remarkable properties of the anticonformal
reflection (abbreviated ACR) w.r.t. $\Ga$.

\medskip\noindent
{\bf Lemma 3.12.} {\em Let $\Omega$ be an R-domain and $R$ the
anticonformal reflection function of its boundary $\Ga$. Then:
\begin{enumerate}
\item[(i)]{If $z\in\wh\C\ssm(\oli\Om\cup B)$, $\{R(z)\}\sbs\Om$}
\smallskip
\item[(ii)]{If $z\in(\Ga\ssm B)$, then one point of $\{R(z)\}$ (namely $z$
    itself) lies on $\Ga$, and the remaining points lie in $\Om$.}
\end{enumerate}}

\medskip

In words, if we let $n$ denote the degree of the rational mapping
$\vap$ (by this we
  mean, writing $\vap=p/q$ where $p$ and $q$ are polynomials without
  common zeroes, that $n$ is the larger of $\deg p$ and $\deg q$),
  then (i) says: all $n$ anticonformal image points of any point outside
$\oli\Om\cup B$ lie in $\Om$. Needless to say, this is a very
special and remarkable property, not shared by all Jordan domains
with algebraic boundary.

\begin{proof}[Proof of Lemma $3.12$] The Schwarz function $S$ of $\Ga$ satisfies
\begin{equation}
S\big(\vap(\ze)\big)=\oli{\vap( \ze)},\qquad |\ze|=1.\tag{3.13}
\end{equation}

Defining
\begin{equation}
\vap^\#(\ze):=\oli{\vap(\bar\ze)},\qquad\ze\in\C\tag{3.14}
\end{equation}
(so that $\vap^\#$ is again a rational function of degree $n$), (3.13)
becomes
\begin{equation}
S\big(\vap(\ze)\big)=\vap^\#(1/\ze)\tag{3.15}
\end{equation}
initially for $|\ze|=1$ and, by analytic continuation globally.
Or, in other words, $S$ is represented parametrically by
\begin{equation}
t=\vap(\ze),\qquad S(t)=\vap^\#(1/\ze).\tag{3.16}
\end{equation}
In terms of the ACR, $R(z):=\oli{S(z)}$, (3.15) becomes
\begin{equation}
R\big(\vap(\ze)\big)=\vap\left (1/\bar \ze\right ),\tag{3.17}
\end{equation}
or parametrically
\begin{equation}
t=\vap(\ze),\qquad R(t)=\vap\left (1/\bar \ze\right).\tag{3.18}
\end{equation}

Let us fix a point $t\in\wh\C\ssm(\oli\Om\cup B)$, and compute
$\{R(t)\}$. Solving $\vap(\ze)=t$ gives $n$ roots
$\ze_1,\dots,\ze_n$ in $\wh\C$ and all $\ze_j$ satisfy
$|\ze_j|>1$, since for $|\ze|\le1$ we have $\vap(\ze)\in\oli\Om$.
Therefore the reflected points $\{1/\bar \zeta_j\colon j=1,\ldots,
n\}$ all lie in the open unit disk ${\mathbb D}$, which implies
that the collection of points $\{\vap(1/\bar \ze_j)\colon j=1, 2,
\dots, n\}$ all lie in ${\Om}$. By (3.18), this proves assertion
(i) of the lemma.

Next, consider $t\in\Ga$. The roots of $\vap(\ze)=t$ consist of
one point, say $\ze_1$, on the unit circle and $n-1$ others
necessarily outside $\oli\D$, call them $\ze_2, \dots, \ze_n$.
Then the elements of $\{R(t)\}$ are $\{\vap(1/\bar\ze_j):j=1, 2,
\dots, n\}$ and the first of these lies on $\Ga$, the remaining
ones in $\Om$. This proves (ii) of the lemma.
\end{proof}

\subsection{Main result} The main result in this section is the
following.

\medskip\noindent
{\bf Theorem 3.19.} {\em If $\Om$ is an $R$-domain of degree $\ge2$,
  the matching problem has no non-trivial solution.}

\begin{proof}

Let $S(z)$ denote the Schwarz function of $\Ga=\partial\Om$ and,
as usual, $R=\oli S$ the ACR associated to $\Ga$. It follows from
(3.16) that the direct analytic continuation of $S$ from $\Ga$ to
$\Om$ (i.e.\ without leaving $\Om$) is meromorphic in $\Om$.
Therefore the analytic continuation of $S$ from $\Ga$ to
$\Om':=\wh\C\ssm\oli\Om$ must encounter a branch point, otherwise
$S$ would be meromorphic in all of $\C$ and this implies, by a
theorem of P.~Davis, that $\Ga$ is a circle (i.e. an $R$-domain of
degree 1, contrary to hypothesis).

Therefore, there is a smooth Jordan arc $\al$ starting from a
point $\ze_0\in\Ga$ and returning to $\ze_0$, with
$\al\ssm\{\ze_0\}\sbs\Om'$ (where $\Omega'$, we recall, denotes
$\wh \C\setminus \oli{\Om}$) such that, by analytic continuation
of $S$ around $\al$ it returns to $\ze_0$ as a different branch,
which locally (near $\ze_0$) we shall denote by $S_1$. By Lemma
3.12 (ii), $\oli{S_1(\ze_0)}\in\Om$, or, expressing matters
henceforth in terms of $R$: the (anti-)analytic continuation of
$R$ from $\ze_0$ to $\ze_0$ along $\al$ leads to a new branch
$R_1$ at $\ze_0$, with $R_1(\ze_0)\in\Om$.

Now, suppose the matching problem has a solution $$f\in
A(\Ga),\qquad g\in A'(\Ga);$$ recall that these functions are
holomorphic in $\Om$, $\Om'$ respectively with $g(\infty)=0$ and
$\oli{f(z)}=g(z)$ for $z\in\Ga$. Thus,
\begin{equation}
\oli{f\big(R(\ze)\big)}=g(\ze)\qquad\ze\in\Ga\tag{3.20}
\end{equation}
where $R$ denotes the ``principal'' branch of the ACR near $\Ga$,
i.e. $R(\ze)=\ze$ on $\Ga$.

If we let $\ze$ describe $\al$ starting from $\ze_0$, (3.20) makes
sense and continues to hold, because $R(\ze)\in\Om$, by Lemma 3.12
(i), where $f$ is defined and holomorphic, and $\oli{f(R(z))}$ is
holomorphic as long as $R(z)$ (which itself is anti-holomorphic)
lies in the domain of holomorphy of $f$.

Now, let $\ze$, after reaching $\ze_0$ along $\al$, continue to
move around $\Ga$ (making small detours into $\Om'$ around any
possible branch points of the branch $R_1$). Then (3.20) continues
to make sense (with $R$ now replaced by $R_1$) and be valid. As
$\ze$ traverses $\Ga$, $R_1(\ze)$ remains ``trapped'' in $\Om$ and
traverses an arc $\be$, starting from $R_1(\ze_0)\in \Om$, which
is a compact subset of $\Om$. This implies, in view of (3.20):
$$
%\underset{\ze\in\Ga}{\max\big|g(\ze)\big|}
\max_{\ze\in\Ga} \big|g(\ze)\big|
 =
\max_{z\in\be} \big|f(z)\big|.
%\underset{z\in\be}{\max\big|f(z)\big|}.
$$
But $|g(\ze)|=|f(\ze)|$ on $\Ga$, and so
$$
\max_{\ze\in\Ga} \big|f(\ze)\big|
=
\max_{\ze\in\be} \big|f(z)\big|.
%\underset{\ze\in\Ga}{\max\big|f(\ze)\big|}
% =
%\underset{z\in\be}{\max\big|f(z)\big|}.
$$
Hence the maximum modulus theorem implies $f$ is constant. This proves
the theorem.
\end{proof}

A similar argument applies to some other domains, for example:

\medskip\noindent
{\bf Theorem 3.21.} {\em If $\Ga$ is a non-trivial ellipse (i.e.\
not a circle), the
  matching problem has only the trivial solution.}

\begin{proof}

The argument is essentially the same as for Theorem 3.19, except
for interchanging the roles played by $\Om$ and $\Om'$. (Indeed,
now $\Om'$ is a ``quadrature domain'', in the sense that it is a
rational conformal map of $\{|\ze|>1\}$, and the Schwarz function
of $\Ga$ is meromorphically extendable to $\Om'$.) We only sketch
the details. $S$ has branch points at the foci of the ellipse,
hence we can continue $R=\oli S$ around a path $\al$ from some
point $\ze\in\Ga$ to $\ze_0$ with $\al\ssm\{\ze_0\}\sbs\Om$, such
that $R$ changes branch upon return to $\ze_0$. Reasoning like
that used earlier shows that $R(z)\in\Om'$ for
$z\in\al\ssm\{\ze_0\}$ and the new branch $R_1$ satisfies
$R_1(\ze_0)\in\Om'$. As $\ze$ now traverses $\Ga$, $R_1(\ze)$
remains ``trapped'' in $\Om'$. (For detailed calculations, see
[E]). A contradiction to the supposition of existence of a
non-trivial matching pair is reached as before.
\end{proof}

\section{Further remarks}

We are aware that our results are extremely limited, and do not
come anywhere near deciding when the ``matching problem'' has
non-trivial solutions, even for domains bounded by algebraic
curves. The following remark pinpoints the delicacy of the
problem.

\medskip

Consider first a polynomial lemniscate: it is defined by
$\{|p(z)|=c\}$ for some $c>0$, where $p$ is a polynomial whose leading
coefficient we may assume is 1. Thus, if
$$
p(z)=(z-z_1) \dots (z-z_n)
$$
the lemniscate is defined by
$$
\prod^n_{j=1} (z-z_j) (\bar z-\bar z_j) - c^2=0
$$
or, writing $z=x+iy$:
$$
P(x,y)=0
$$
where
$$
P(x,y)= \prod^n_{j=1} (x^2+y^2+a_j x+b_jy + c_j) - c^2,
$$
and $a_j$, $b_j$, $c_j$ are real constants (depending on $z_j$).

Hence,
\begin{equation}
P(x,y) = (x^2+y^2)^n + \dots \,\,
 \textrm{(terms of degree $<2n$)}.
\end{equation}
Now, on the other hand, it is easy to check that {\em every
$R$-domain
  is bounded by a curve $\{P=0\}$ where $P$ has the form $(4.1)$ for
  some $n$}.

Thus, among domains bounded by curves $\{P=0\}$ where $P$ has the form
(4.1) for some $n\ge2$, there exist domains for which the matching
problem admits non-trivial solutions, as well as ones for which this is
not the case. This shows that the issue cannot be decided just by
examining the highest degree terms that appear in $P$.

\medskip

Our next remark concerns a result that follows by comparison of
Theorems 2.1 and 3.19. We formulate it as:

\medskip\noindent
{\bf Theorem 4.2.} {\em Let $P$ be a monic polynomial, and suppose the
  (lemniscate) domain $\Om:=\{z:\ |P(z)|<1\}$ is a Jordan domain, which
  also is a quadrature domain. Then $\Om$ is a disk.}

\medskip

It seems of some interest to give a simple self-contained proof of
this result, which we proceed to do:

\medskip

Suppose the quadrature domain $\Om$ has order $m$, so $S$, the
Schwarz function of its boundary, has a branch equal to $\bar z$
on $\Ga=\po$, and extendible meromorphically throughout $\Om$ with
$m$ poles (for background, see [S]). Now, on $\Ga$, $|P(z)|^2=1$,
i.e. $1=P(z)\oli{P(z)}=P(z)P^\#(S(z))$. By analytic continuation,
this holds for all $z$ in $\Om$. This implies $P^\#(S(z))$ has
$k=\deg P$ poles in $\Om$. But $P^\#(S(z))$ has $k m$ poles
(counting always with multiplicities), so $m=1$, which implies
that $\Om$ is a disk. This fact is well known, but for the
convenience of the reader we give here the argument. If $m=1$ then
$S(z)$ has a simple pole at, say, $a\in \Om$. Thus, the function
$Q(z):=(z-a)(S(z)-\bar a)$ has a holomorphic extension to a
neighborhood of $\overline{\Om}$. Moreover, on $\Ga$ we have
$Q(z)=(z-a)(\bar z-\bar a)=|z-a|^2$ and, in particular, $Q(z)$ is
real-valued on $\Ga$ . It follows that $Q(z)$ is constant in $\Om$
and, hence, $|z-a|^2$ is constant on $\Ga$, which implies that
$\Om$ is a disk.

\medskip

In conclusion, let us briefly take note of some natural extensions of
the ``matching problem'':

\medskip\noindent
(i) {\em Extension to $n>2$ variables.} The operator we called
    $\Pi$ is, in essence, representable as the two-dimensional case
    of the operator $I+2K$, where $K$ is the integral operator mapping a function
    $\varphi$ on $\Ga$ to $u$,
    where
$$
u(y):=\il_\Ga \varphi(x) (\partial/\partial N_x) E(y-x)
dS_x,\qquad y\in\Ga
$$
and $E$ is the Newtonian kernel. Our ``matching problem'' was
equivalent to whether this operator (when $n=2$) is {\em
injective} (in one or another space of functions on $\Ga$). This
problem makes good sense in $\R^n$ for $n\ge 3$ also, but we have
no results except for the easily verified fact that for the sphere
in $\R^n$, with $n>2$, the operator is indeed injective unlike the
case $n=2$. For in higher dimensions the kernel of the operator
for a sphere is a constant times the Newtonian kernel, so a
nontrivial nullspace would imply existence of a function whose
single layer potential vanishes on the sphere and hence by the
maximum principle throughout the space, an obvious contradiction.

To see that the kernel $k(x,y):=(\partial/\partial N_x) E(y-x)$ of
the operator $K$ is equal to a constant times $E(y-x)$ when $n\geq
3$ and $\Ga$ is the unit sphere, we recall that
\begin{equation}
E(y-x)=c_n\frac{1}{|y-x|^{n-2}}
\end{equation}
and so, when $\Ga$ is the unit sphere in $\R^n$,
\begin{equation}
k(x,y)=-(n-2)c_n\frac{(y-x)\cdot x}{|y-x|^{n}}.
\end{equation}
Now, for $x,y$ on the unit sphere, we have
\begin{equation}
\begin{aligned}
|y-x|^2 &=(y-x)\cdot (y-x)\\
&= 1-x\cdot y-y\cdot x+1\\
&=2(1-x\cdot y)\\
&=2(x\cdot x-x\cdot y)\\
&= -2(y-x)\cdot x
\end{aligned}
\end{equation}
and, hence,
\begin{equation}
k(x,y)=-\frac{n-2}{2} E(y-x),
\end{equation}
as claimed.

\smallskip\noindent
(ii) {\em Generalized matching.} In the case $n=2$ consider the
    following problem: $\Om$ is a given Jordan domain with smooth
    boundary $\Ga$, and we seek two functions $f_1$, $f_2$ holomorphic
    respectively in $\Om$ and $\Om':=\wh\C\ssm\oli\Om$, with
    $f_2(\infty)=0$ and such that the following relations hold on
    $\Ga$, where $f_j=u_j+iv_j$ $(j=1, 2)$ with $u_j$, $v_j$
    real-valued:
\begin{equation*}
\begin{matrix}
A_1 v_1 &+ &B_1 v_1 &+ &A_2 u_1 &+ &B_2 v_2 &= \vap\\
C_1 v_1 &+ &D_1 v_1 &+ &C_2 u_2 &+ &D_2 v_2 &= \psi
\end{matrix}
\tag{$*$}
\end{equation*}
where $A_j, \dots, D_j$ ($j=1, 2$) are given real constants and
$\vap$, $\psi$ given real functions on $\Ga$. Observe that our
``matching problem'' is the special case where this system reduces
to $u_1-u_2=0$, $v_1+v_2=0$. For varying choices of the
parameters, the system ($*$) becomes the model for a wide range of
problems in accoustics, electrostatics, gravitational attraction,
diffraction theory, etc., see [C] for details and references.
Still further generalizations are to allow $A_j$ etc. to be
functions on $\Ga$. Also, we can study multivariable analogs where
$u_1$, $v_2$ are replaced by $U$, $\partial U/\partial N$ ($U$
being harmonic in $\Om$) and $u_2$, $v_2$ by $V$, $\partial
V/\partial N$ where $V$ is harmonic in the exterior domain with
$V(\infty)=0$. Few of these problems are completely solved.

\vfill

\end{document}